# The Number of Subgroups and Cyclic Subgroups of Finite Group and Its Application by GAP Program

Abdallah Shihadeh

Department of Mathematics, Faculty of Science, The Hashemite University, Zarqa 13133, PO box 330127, Jordan

abdallaha_ka@hu.edu.jo

.

**Abstract.**

In this work, In this paper, we presented A novel approach for calculating a finite group's set of subgroups, cyclic subgroups using it to establish the quantity of all subgroups in the direct product of two groups, the Dicyclic group have order $4n$ and the Cyclic group have order $p$, $\tau(n)$ is a total number of all divisors of $n$ and the summation of all divisors, $\sigma(n)$. Using this, we derive a formula for the number of subgroups in the group $T_{4n} \times C_p$. We then use the computer program GAP to find all $T_{4n} \times C_p$ with exactly $|T_{4n} \times C_p|$ – t cyclic subgroups for $t \geq 1$.

## 1- Introduction

CGT, or Computational Group Theory, is a well-established branch of algebra in the field of computation. While many symbolic algebra programs can work with groups to some degree, CGT specifically studies both theoretical and practical algorithms for dealing with groups. These algorithms have applications not only in group theory, but in other fields such as combinatorics, topology, and physics, where groups are used to represent symmetries. For these computations, two systems are particularly well-suited: GAP and Magma. Additionally, there are many other independent programs and smaller systems that are also accessible.

In this paper for all groups are finite. Let's suppose we have a group $G$ of finite order $n$ We want to algorithmically compute all subgroups, and there are some ways to do that.

In their respective presentations, Cavior and Calhoun [1] and [2] showed the number of all subgroups in the one of finite group is known the dihedral group to be $\tau(n) + \sigma(n)$. Ahrafi with authors introduced a new rule to calculate the number of all subgroups for same of a finite groups, citing earlier research [12,13], this method depend on divide of the number of $n = 2^r m$, where $m = \prod_{i=1}^{r} p_i^{\zeta_i}$ to integer divisors and structure the table have first row on all divisors of $2^r$ and first column have all divisors of $m = \prod_{i=1}^{r} p_i^{\zeta_i}$ to find the order of subgroups.

| * \ K | 1 | 2 | ... | $r$ | $r+1$ |
|---|---|---|---|---|---|
|  | 1 | 2 | ... | $2^{r-1}$ | $2^r$ |
| 1 | 1 | $2*1$ | ... | $2^{r-1}*1$ | $2^r*1$ |
| $p$ | $p$ | $2*p$ | ... | $2^{r-1}*p$ | $2^r*p$ |
| $\vdots$ | $\vdots$ | $\vdots$ | ... | $\vdots$ | $\vdots$ |
| $m$ | $m$ | $2*m$ | ... | $2^{r-1}*m$ | $2^r*m$ |

Lazorec with M. Trume [4,14] and M. S. [15,16,17] emphasized the importance of the cyclic subgroup count in establishing the total quantity of a group in [7]. They presented that for a finite group G with its poset of cyclic subgroups $L_1(G)$ [1,3,6,7] the quantity is $\propto (G) = \frac{|L_1(G)|}{|G|}$ Within this framework, a finite group is considered nilpotent if $\propto (G)$ is equal to $3/4$. [2,5,8] Building upon this research, this paper aims to delve deeper into determining the number of subgroups in the direct product of two groups, specifically $\Phi \simeq T_{4n} \times C_p$, where $p$ is an odd of all prime number and $p$ is not a divisor of $n$. G. James. M. Liebeck in [7] presented the dicyclic group by $T_{4n} = \langle \alpha, \beta | \alpha^{2n} = e, \alpha^n = \beta^2 = e, \beta\alpha\beta^{-1} = \alpha^{-1} \rangle$ has order $4n$.

This paper is divided into sections for easy understanding. Initially, we elaborate on the core concept of displaying subgroups of finite groups.

Following that our paper introduces a crucial algorithm: a technique for exhibiting specific subgroups in a finite group and we demonstrate how to utilize these subgroups to compute head complements, while, we explain their application in computing formation-theoretic subgroups. To highlight the practicality of our methods, we present the implementation results of these techniques using the computer algebra system GAP (see Schonert et al., 1995). In [9,10,12,15] the final sectionwe provide a comprehensive analysis of the complexity involved. A cyclic group $G$ is a group that can be produced by just one element $\alpha$, resulting in every element in $G$ taking the form $\alpha^i$ with an integer $i$. $C_n$ represents the cyclic group of order $n$.

We frequently encounter cyclic groupings in daily life. A group that has an element to which an operation is done and which yields the entire set is called a *c*yclic group. The *s*implest *g*roup is a *c*yclic *g*roup. A natural pattern or a geometric design we create ourselves can both be examples of cyclic groups. In [4,6,11,16] Cyclic groups can also be thought of as reels since, given enough rotations, each item will eventually return to its initial location. In this research study, we investigate more cyclic group operations in number proposition, such as the division algorithm and Chinese remainder theorem, as well as other operations such as bell ringing, modular systems, chaotic proposal, 12-Hours wristwatch, and direct codes. But they might utilize the inventor to identify the quickest straightforward circuit for. For more information, see [5,7,9,11,12].

**Theorem[6]:**

A cyclic group is cyclic in all of its subgroups. If $G = \langle \alpha \rangle$ is cyclic, then $\alpha^{|G|/d}$ can generate by only one of the subgroup have order $d$ for each divisor $d$ of $|G|$.

**Theorem[6]:**

Every group of composite order has proper subgroups.

2- **Main Results**

The our idea in the section, we will provide a novel approach to determine how many subgroups there are in any order division order group. $T_{4n} \times C_p$. This algorithm will help in calculating these subgroups through the number of the subgroups table proposed in the following results, which will quickly reach the result while determining the order of these subgroups.

We will in this section reviewer some of the standard of claims using information from [1, 2], Cavior was able to get the bulding of the $d$ihedral group's subgroups. and shelash in 2020 computed the set of all subgroups of the dihedral group with some of cyclic groups.

**Lemma 2.1.**

Suppose that $T_{4n}$ be a $D$icyclic group have $4n$ order and and $C_p$ is a $C$yclic finite groups, the presentation of the direct product of two groups is given by:

$$\Phi \cong T_{4n} \times C_p = \langle \alpha, \beta, c | \alpha^{2n} = e, \alpha^n = \beta^2 = \gamma^p = e, \beta\alpha\beta^{-1} = \alpha^{-1}, [\alpha,\gamma] = [\beta,\gamma] = 1 \rangle$$

**Proof:**

It is clear that, the Dicyclic group $T_{4n}$ is define by $T_{4n} = \langle \alpha, \beta | \alpha^{2n} = e, \alpha^n = \beta^2 = e, \beta\alpha\beta^{-1} = \alpha^{-1} \rangle$, it is generated by $\alpha$ and $\beta$ are elements, where order $O(\alpha) = 2n$ and the order $O(\beta) = 4$. So, the cyclic group is define by $C_n = \langle \gamma | \gamma^n = e \rangle$, By defining the direct product of two groups we can obtain the following.

$$T_{4n} \times C_p = \langle \alpha, \beta | \alpha^{2n} = e, \alpha^n = \beta^2 = e, \beta\alpha\beta^{-1} = \alpha^{-1} \rangle \times \langle \gamma | \gamma^p = e \rangle$$

$$= \langle \alpha, \beta, \gamma | \alpha^{2n} = e, \alpha^n = \beta^2 = \gamma^p = e, \beta\alpha\beta^{-1} = \alpha^{-1}, \alpha\gamma = \gamma\alpha, \beta\gamma = \gamma\beta \rangle$$

$$= \langle \alpha, \beta, \gamma | \alpha^{2n} = e, \alpha^n = \beta^2 = \gamma^p = e, \beta\alpha\beta^{-1} = \alpha^{-1}, [\alpha,\gamma] = [\beta,\gamma] = e \rangle$$

Hence $[x, y] = xyx^{-1}y^{-1}$ for each $x, y \in G$.

Thus, $\alpha\gamma = \gamma\alpha \Rightarrow \alpha\gamma\alpha^{-1}\gamma^{-1} = e \Rightarrow [\alpha,\gamma]$, so for $\beta\gamma = \gamma\beta \Rightarrow \beta\gamma\beta^{-1}\gamma^{-1} = e \Rightarrow [\beta,\gamma]$.

**Proposition 2.2.**

Suppose that $n = 2^r \prod_{i \in \mathbb{N}} p_i^{\zeta_i}$ be a positive integer number. Then:

1) $\mathcal{H}_1 = \langle \alpha^i \rangle \ \forall \ i|2n$;
2) $\mathcal{H}_2 = \langle \alpha^i, \alpha^j \beta \rangle$, $\forall \ i|n$ and $1 \leq j \leq i$;
3) $\mathcal{H}_3 = \langle \alpha^i, \gamma \rangle$, $\forall \ i|2n$;
4) $\mathcal{H}_4 = \langle \alpha^i, \alpha^j \beta, \gamma \rangle$, $\forall \ i|n$ and $1 \leq j \leq i$.

Are subgroup of $T_{4n} \times C_p$.

***Proof.***

From [6,7,8,9] the set of all subgroups of the dicyclic group it is given by the strucuter presented by $\mathcal{H}_1$ and $\mathcal{H}_2$.

Now, we will prove the structure subgroups for $\mathcal{H}_3$ and $\mathcal{H}_4$. Becouse the cyclic group $C_p$ is a prime order then it is has cases, in case one is $\langle e \rangle$, case two is $\langle \gamma \rangle$, for case $\langle e \rangle$, we obtaian on the following:

$$\mathcal{H}_1 \times \langle e \rangle = \mathcal{H}_1 \text{ for each } i|2n$$

$$\mathcal{H}_2 \times \langle e \rangle = \mathcal{H}_2 \text{ for each } i|n$$

And for case two, suppose that $\alpha$ and $\gamma$ are elements in abelian groups $\mathcal{H}_1$ and $\mathcal{H}_3$ respectively also, their orders are limited as well as co-prime, then $\mathcal{H}_1 \times \mathcal{H}_3 = \langle \alpha, \gamma \rangle$
We have the order $((\alpha, \gamma)) = lcm\{o(\alpha), o(\gamma)\} = o(\alpha)o(\gamma)$, Consequently subgroup $(\langle \alpha, \gamma \rangle)$ has order $o(\alpha)o(\gamma)$, which same request of $\langle \alpha \rangle \times \langle \gamma \rangle$. Presently need to demonstrate that $\langle \alpha, \gamma \rangle = \langle \alpha \rangle \times \langle \gamma \rangle$ Suppose that $\mu \in (\alpha, \gamma)^\omega$ be any inconsistent component of $\langle \alpha, \gamma \rangle$, $\mu = (\alpha^\omega, \gamma^\omega) \Longrightarrow \mu \in \alpha^\omega$, and $\mu \in \gamma^\omega$ thus $(\alpha^\omega, \gamma^\omega) \in \langle \alpha \rangle \times \langle \gamma \rangle$ and in final $\mu \in \langle \alpha \rangle \times \langle \gamma \rangle$.

Hence $\langle \alpha, \gamma \rangle \subseteq \langle \alpha \rangle \times \langle \gamma \rangle$. So $\langle \alpha \rangle \times \langle \gamma \rangle \subseteq \langle \alpha, \gamma \rangle$ by similarly, we can obtain on the $\langle \alpha, \gamma \rangle = \langle \alpha \rangle \times \langle \gamma \rangle$ thus for each $i|2n$ there exist a cyclic subgroup of type $\langle \alpha^i, \gamma \rangle$.

$$\mathcal{H}_1 \times \langle \gamma \rangle = \langle \alpha^i, \gamma \rangle = \mathcal{H}_3 \text{ for each } i|2n$$

$$\mathcal{H}_2 \times \langle \gamma \rangle = \langle \alpha^i, \alpha^j\beta, \gamma \rangle = \mathcal{H}_4 \text{ for each } i|n$$

## 3- Novel Formula in Number Theory for the quantity of subgroups of $T_{4n} \times C_p$ where $p \nmid n$

**Theorem 3.1** (Order Subgroup Table)

Suppose that $n = 2^r m$ is an integer number and $m = p_1^{\zeta_1} p_2^{\zeta_2} \dots p_j^{\zeta_j}$ is the prime factorization of $j > 1$. The number subgroups order is given by the following table:

1- Fisrt Part of Table: if $\lambda|m$, then:

| | k | 1 | 2 | 3 | | r+2 | r+3 |
|---|---|---|---|---|---|---|---|
| | * | 1 | 2 | $2^2$ | ... | $2^{r+1}$ | $2^{r+2}$ |
| | 1 | 1 | 1 | $n+1$ | ... | $\frac{n}{2^{r-1}}+1$ | $\frac{n}{2^r}$ |
| $\lambda\|m$ | ⋮ | ⋮ | ⋮ | ⋮ | ... | ⋮ | ⋮ |
| | m | 1 | 1 | $\frac{n}{m}+1$ | ... | $\frac{n}{2^{r-1}m}+1$ | $\frac{n}{2^r m}$ |
| Sum First part | | $\tau(m)$ | $\tau(m)$ | $2^r \sigma(m) + \tau(m)$ | ... | $2\sigma(m) + \tau(m)$ | $\sigma(m)$ |

2- Fisrt Part of Table: if $p(\lambda|m)$, then:

| | k | 1 | 2 | 3 | | r+2 | r+3 |
|---|---|---|---|---|---|---|---|
| | * | 1 | 2 | $2^2$ | ... | $2^{r+1}$ | $2^{r+2}$ |
| | p | 1 | 1 | $(\frac{np}{p}+1)$ | ... | $(\frac{np}{2^{r-1}p}+1)$ | $(\frac{np}{2^r p})$ |
| $p(\lambda\|m)$ | ⋮ | ⋮ | ⋮ | ⋮ | ... | ⋮ | ⋮ |
| | mp | 1 | 1 | $(\frac{np}{mp}+1)$ | ... | $(\frac{np}{2^{r-1}mp}+1)$ | $(\frac{np}{2^r mp})$ |

| Sum Second Part | $\tau(m)$ | $\tau(m)$ | $2^r\sigma(m) + \tau(m)$ | ... | $2\sigma(m) + \tau(m)$ | $\sigma(m)$ |

*Proof.*

In the first, we will calculate the total number of subgroups for columns with $k = 1,2$. This is equivalent to the number of positive divisors of a positive integer $m$, denoted as $\tau(m)$. Let $t$ and $m$ be positive integers, where $t$ is a divisor of $m$. As $t$ is a divisor of $m$, there exists an integer $l$ such that $m/l = t$. The set $\{\rho^{\wedge}t, \rho^{\wedge}2t, \rho^{\wedge}3t, \dots, \rho^{\wedge}(m-t), \rho^{\wedge}m\}$ will be formed when $\rho^{\wedge}t$ is raised to various powers. For this set to be a subgroup, each power of $\rho$ must generate a subgroup. This means that every divisor of m will generate a subgroup of rotations. Furthermore, any multiple of $t$ that is not a divisor of m will generate the same subgroup as $\rho^{\wedge}t$, and any power of $\rho$ that is relatively prime to $m$ will also generate the same subgroup as $\rho$. In conclusion, the number of subgroups for $k = 1,2$ will be $\tau(m)$ for both cases.

Now, for $k = 3, 4, \dots, r+2$, we can explain this by the $\frac{n}{2^{k-1}t} + 1$ for each $t|m$ are computing by the following:

$$\frac{n}{2^{k-3}} + 1 + \frac{n}{2^{k-3}p} + 1 + \cdots + \frac{n}{2^{k-3}m} + 1$$

$$= 2^{r-k+3}\left(n + \frac{n}{p} + \cdots + 1\right) + (1 + 1 + \cdots + 1)$$

$$= 2^{r-k+3}\sigma(m) + \tau(m) \quad \text{for each} \quad 3 \leq k \leq r+2 \quad \text{and for last column are equal to}$$

$\sigma(m)$.

Now, the summation set of the subgroups of all subgroups for any columns is given by:

$$\tau(m) + \tau(m)2^r\sigma(m) + \tau(m) + 2^{r-1}\sigma(m) + \tau(m) \ldots + 2\sigma(m) + \tau(m) + \sigma(m) =$$

$$\sigma(2^r)\sigma(m) + r\tau(m) = \sigma(n) + (r+2)\tau(m) = \sigma(n) + \tau(2n).$$

**Theorem 3.2**

Let $n = 2^r \prod_{i \in \mathbb{N}} p_i^{\zeta_i}$ be a positive integer number. The number of subgroups of group $\Phi$ is given by:

$$\Sigma NS(\Phi) = \tau(p)(\tau(2n) + \sigma(n))$$

*Proof.*

From **Theorem (2.2),** we can computing the all of the subgroups from using the number subgroups of table, the summation of columns of the first part we obtain on the following :

$Sub(Fir.Part)$

$$= \tau(m) + 2^r\sigma(m) + \tau(m) + \cdots + 2^2\sigma(m) + \tau(m) + 2\sigma(m)$$
$$+ \tau(m) + \sigma(m)$$

$$= \underbrace{\tau(m) + \cdots + \tau(m)}_{r+2} + 2^r\sigma(m) + \cdots + 2\sigma(m) + \sigma(m)$$

$$= \tau(2^{r+1})\tau(m) + (2^r + 2^{r-1} + \cdots + 2 + 1)\sigma(m)$$

$$= \tau(2^{r+1})\tau(m) + \sigma(2^r)\sigma(m)$$

$$= \tau(2n) + \sigma(n).$$

Comparably, the total of all subgroups in the second section of Table 2.1 should equal $\tau(n) + \sigma(n)$. Since the first table contains all order subgroups that $p$ does not divide and the second contains all order subgroups that $p$ divides, the number parts of Table 2 should equal $\tau(p) = 2$.

In the final we obtain on the

$$\Sigma NS(\Phi) = 2(\tau(2n) + \sigma(n)) = \tau(p)(\tau(2n) + \sigma(n))$$

**Example 3.3**

In the following example, we will present an application for computing the number of subgroups order table, using the Theorem [3.1]

Take $n = 23^2 5^2$ and $p = 7$, thus, the direct product group $\Phi = T_{1800} \times C_7$ thave order 12600, the number of all subgroups are given by the following table:

1) $k = 1,2,3,4$
2) The set of divisors of $n$ are {1,2,3,5,6,9,10,15,18,25,30,45,50,75,90,150,225,450}
3) Order of subgroups table:
4) If ג|225, then :

|  | K | 1 | 2 | 3 | 4 |
|---|---|---|---|---|---|
|  | * | 1 | 2 | $2^2$ | $2^3$ |
| ג|225 | 1 | 1 | 1 | 451 | 225 |
|  | 3 | 1 | 1 | 151 | 75 |
|  | 5 | 1 | 1 | 91 | 45 |
|  | 9 | 1 | 1 | 51 | 25 |
|  | 15 | 1 | 1 | 31 | 15 |
|  | 25 | 1 | 1 | 19 | 9 |
|  | 45 | 1 | 1 | 11 | 5 |

|  | 75 | 1 | 1 | 7 | 3 |
|---|---|---|---|---|---|
|  | 225 | 1 | 1 | 3 | 1 |
| Sum First part | | 9 | 9 | $2(403) + 9$ | 403 |

5) If $7(\lambda|225)$, then:

|  | K | 1 | 2 | 3 | 4 |
|---|---|---|---|---|---|
|  | * | 1 | 2 | $2^2$ | $2^3$ |
| $7(\lambda|225)$ | 7 | 1 | 1 | 451 | 225 |
|  | 21 | 1 | 1 | 151 | 75 |
|  | 35 | 1 | 1 | 91 | 45 |
|  | 63 | 1 | 1 | 51 | 25 |
|  | 105 | 1 | 1 | 31 | 15 |
|  | 175 | 1 | 1 | 19 | 9 |
|  | 315 | 1 | 1 | 11 | 5 |
|  | 525 | 1 | 1 | 7 | 3 |
|  | 1575 | 1 | 1 | 3 | 1 |
| Sum Second Part | | 9 | 9 | $2(403) + 9$ | 403 |

To clarify how these subgroups are calculated. We will calculate all subgroups based on the following Proposition [2.2].

From this type subgroup $\mathcal{H}_1 = \langle \alpha^i \rangle \ \forall \ i|900$ we get the following subgroups:

| Subgroups | Order | Subgroups | Order | Subgroups | Order |
|---|---|---|---|---|---|
| $\langle \alpha \rangle$ | 900 | $\langle \alpha^{15} \rangle$ | 60 | $\langle \alpha^{75} \rangle$ | 12 |
| $\langle \alpha^2 \rangle$ | 450 | $\langle \alpha^{18} \rangle$ | 50 | $\langle \alpha^{90} \rangle$ | 10 |

| Subgroup | Order | Subgroup | Order | Subgroup | Order |
|---|---|---|---|---|---|
| $\langle \alpha^3 \rangle$ | 300 | $\langle \alpha^{20} \rangle$ | 45 | $\langle \alpha^{100} \rangle$ | 9 |
| $\langle \alpha^4 \rangle$ | 225 | $\langle \alpha^{25} \rangle$ | 36 | $\langle \alpha^{150} \rangle$ | 6 |
| $\langle \alpha^5 \rangle$ | 180 | $\langle \alpha^{30} \rangle$ | 30 | $\langle \alpha^{180} \rangle$ | 5 |
| $\langle \alpha^6 \rangle$ | 150 | $\langle \alpha^{36} \rangle$ | 25 | $\langle \alpha^{225} \rangle$ | 4 |
| $\langle \alpha^9 \rangle$ | 100 | $\langle \alpha^{45} \rangle$ | 20 | $\langle \alpha^{300} \rangle$ | 3 |
| $\langle \alpha^{10} \rangle$ | 90 | $\langle \alpha^{50} \rangle$ | 18 | $\langle \alpha^{450} \rangle$ | 2 |
| $\langle \alpha^{12} \rangle$ | 75 | $\langle \alpha^{60} \rangle$ | 15 | $\langle e \rangle$ | 1 |

From this type subgroup $\langle \alpha^i, \alpha^j \beta \rangle$, $\forall\, i|450$ and $1 \leq j \leq i$;

we get the following subgroups:

| Subgroups | Order | Subgroups | Order |
|---|---|---|---|
| $\langle \alpha, \beta \rangle$ | 1800 | $\langle \alpha^{25}, \alpha^j \beta \rangle_{1 \leq j \leq 25}$ | 72 |
| $\langle \alpha^2, \alpha^j \beta \rangle_{1 \leq j \leq 2}$ | 900 | $\langle \alpha^{30}, \alpha^j \beta \rangle_{1 \leq j \leq 30}$ | 60 |
| $\langle \alpha^3, \alpha^j \beta \rangle_{1 \leq j \leq 3}$ | 600 | $\langle \alpha^{45}, \alpha^j \beta \rangle_{1 \leq j \leq 45}$ | 40 |
| $\langle \alpha^5, \alpha^j \beta \rangle_{1 \leq j \leq 5}$ | 360 | $\langle \alpha^{50}, \alpha^j \beta \rangle_{1 \leq j \leq 50}$ | 36 |
| $\langle \alpha^6, \alpha^j \beta \rangle_{1 \leq j \leq 6}$ | 300 | $\langle \alpha^{75}, \alpha^j \beta \rangle_{1 \leq j \leq 75}$ | 24 |
| $\langle \alpha^9, \alpha^j \beta \rangle_{1 \leq j \leq 9}$ | 200 | $\langle \alpha^{90}, \alpha^j \beta \rangle_{1 \leq j \leq 90}$ | 20 |
| $\langle \alpha^{10}, \alpha^j \beta \rangle_{1 \leq j \leq 10}$ | 180 | $\langle \alpha^{150}, \alpha^j \beta \rangle_{1 \leq j \leq 150}$ | 12 |
| $\langle \alpha^{15}, \alpha^j \beta \rangle_{1 \leq j \leq 15}$ | 120 | $\langle \alpha^{225}, \alpha^j \beta \rangle_{1 \leq j \leq 225}$ | 8 |
| $\langle \alpha^{18}, \alpha^j \beta \rangle_{1 \leq j \leq 18}$ | 100 | $\langle \alpha^{450}, \alpha^j \beta \rangle_{1 \leq j \leq 450}$ | 4 |

From this type subgroup $\mathcal{H}_3 = \langle \alpha^i, \gamma \rangle$ $\forall\, i|900$ we get the following subgroups:

| Subgroups | Order | Subgroups | Order | Subgroups | Order |
|---|---|---|---|---|---|
| $\langle \alpha, \gamma \rangle$ | 6300 | $\langle \alpha^{15}, \gamma \rangle$ | 420 | $\langle \alpha^{75}, \gamma \rangle$ | 84 |
| $\langle \alpha^2, \gamma \rangle$ | 3150 | $\langle \alpha^{18}, \gamma \rangle$ | 350 | $\langle \alpha^{90}, \gamma \rangle$ | 70 |
| $\langle \alpha^3, \gamma \rangle$ | 2100 | $\langle \alpha^{20}, \gamma \rangle$ | 315 | $\langle \alpha^{100}, \gamma \rangle$ | 63 |
| $\langle \alpha^4, \gamma \rangle$ | 1575 | $\langle \alpha^{25}, \gamma \rangle$ | 252 | $\langle \alpha^{150}, \gamma \rangle$ | 42 |
| $\langle \alpha^5, \gamma \rangle$ | 1260 | $\langle \alpha^{30}, \gamma \rangle$ | 210 | $\langle \alpha^{180}, \gamma \rangle$ | 35 |
| $\langle \alpha^6, \gamma \rangle$ | 1050 | $\langle \alpha^{36}, \gamma \rangle$ | 175 | $\langle \alpha^{225}, \gamma \rangle$ | 28 |
| $\langle \alpha^9, \gamma \rangle$ | 700 | $\langle \alpha^{45}, \gamma \rangle$ | 140 | $\langle \alpha^{300}, \gamma \rangle$ | 21 |
| $\langle \alpha^{10}, \gamma \rangle$ | 630 | $\langle \alpha^{50}, \gamma \rangle$ | 126 | $\langle \alpha^{450}, \gamma \rangle$ | 14 |
| $\langle \alpha^{12}, \gamma \rangle$ | 525 | $\langle \alpha^{60}, \gamma \rangle$ | 105 | $\langle \gamma \rangle$ | 7 |

From this type subgroup $\langle \alpha^i, \alpha^j \beta, \gamma \rangle$, $\forall\, i | 450$ and $1 \leq j \leq i$;

we get the following subgroups:

| Subgroups | Order | Subgroups | Order |
|---|---|---|---|
| $\langle \alpha, \beta, \gamma \rangle$ | 12600 | $\langle \alpha^{25}, \alpha^j \beta, \gamma \rangle_{1 \leq j \leq 25}$ | 504 |
| $\langle \alpha^2, \alpha^j \beta, \gamma \rangle_{1 \leq j \leq 2}$ | 6300 | $\langle \alpha^{30}, \alpha^j \beta, \gamma \rangle_{1 \leq j \leq 30}$ | 420 |
| $\langle \alpha^3, \alpha^j \beta, \gamma \rangle_{1 \leq j \leq 3}$ | 4200 | $\langle \alpha^{45}, \alpha^j \beta, \gamma \rangle_{1 \leq j \leq 45}$ | 280 |
| $\langle \alpha^5, \alpha^j \beta, \gamma \rangle_{1 \leq j \leq 5}$ | 2520 | $\langle \alpha^{50}, \alpha^j \beta, \gamma \rangle_{1 \leq j \leq 50}$ | 252 |
| $\langle \alpha^6, \alpha^j \beta, \gamma \rangle_{1 \leq j \leq 6}$ | 2100 | $\langle \alpha^{75}, \alpha^j \beta, \gamma \rangle_{1 \leq j \leq 75}$ | 168 |
| $\langle \alpha^9, \alpha^j \beta, \gamma \rangle_{1 \leq j \leq 9}$ | 1400 | $\langle \alpha^{90}, \alpha^j \beta, \gamma \rangle_{1 \leq j \leq 90}$ | 140 |
| $\langle \alpha^{10}, \alpha^j \beta, \gamma \rangle_{1 \leq j \leq 10}$ | 1260 | $\langle \alpha^{150}, \alpha^j \beta, \gamma \rangle_{1 \leq j \leq 150}$ | 84 |
| $\langle \alpha^{15}, \alpha^j \beta, \gamma \rangle_{1 \leq j \leq 15}$ | 840 | $\langle \alpha^{225}, \alpha^j \beta, \gamma \rangle_{1 \leq j \leq 225}$ | 56 |
| $\langle \alpha^{18}, \alpha^j \beta, \gamma \rangle_{1 \leq j \leq 18}$ | 700 | $\langle \alpha^{450}, \alpha^j \beta, \gamma \rangle_{1 \leq j \leq 450}$ | 28 |

Sum first part = 1236 subgroups

Sum second part = 1236 subgroups

So, we found out that there are a total of 1236+1236=2472 subgroups in $\Phi = T_{1800} \times C_7$.

By using the final formula proposed in the Theorem 3.2 to verify the final result, we notice that the number be equal to the number of calculation in the table.

$$\Sigma NS(T_{1800} \times C_7) = \tau(7)(\tau(900) + \sigma(450))$$

$$= (2)(27 + 1209)$$

$$= 2472$$

By GAP program.

```
gap> n:=450;
gap>F := FreeGroup( "a", "b" );
gap>T := F / [ a^(2*n), b^2/a^n, b^(-1)*a*b*a ];
gap> i:=IsomorphismPermGroup(T);
gap>c:=CyclicGroup(7);
gap> Tc:=DirectProduct(T,c);
gap> s:=AllSubgroups(Tc);
gap>Size(s);
Answer:  2472
```

### 4- The number of CyclicGroup of group $T_{4n} \times C_p$ where $p \nmid n$

In this section we will present, the cyclic group is one of important application group theory.

**Proposition 4.1**

A Cyclic subgroups $\langle a^i \rangle$, where $i|2n$; they are isomorphic with $C_{\frac{2n}{i}}$

**Proof:**

Clear that the sturucture subgroup $\langle \alpha^i \rangle$ generated by one element for all $i|2n$ is isomorphic to cyclic subgroup of group $\Phi$. Assume that, there exist $k$ is an integer number divided $2n$, then $(\alpha^k)^{2n/k} = e$, thus mean $\langle \alpha^k \rangle = \{ a^k, a^{2k}, \dots, a^{2n/k} \} \cong C_{\frac{2n}{i}}$,

**Proposition 4.2**

A Cyclic subgroups $\langle \alpha^k, \alpha^j \beta \rangle$, where $k = n$ and $1 \leq j \leq k$ be isomorphic to $C_4$.

**Proof:**

Since $\alpha^{2n} = e$, this mean $\alpha^n$ is order 2, and $\alpha^j \beta$ of order 4, thus the strucutre subgroup $\langle \alpha^n, \alpha^j \beta \rangle = \{e, \alpha^n, \alpha^j \beta, \alpha^{j+n} \beta\}$, be isomorphic to $C_4$ for all $1 \leq j \leq k$

$$\langle \alpha^n, \alpha^1 \beta \rangle = \{e, \alpha^n, \alpha^1 \beta, \alpha^{1+n} \beta\},$$

$$\langle \alpha^n, \alpha^2\beta \rangle = \{e, \alpha^n, \alpha^2\beta, \alpha^{2+n}\beta\},$$

$$\vdots$$

$$\langle \alpha^n, \alpha^n\beta \rangle = \{e, \alpha^n, \alpha^n\beta, \alpha^{n+n}\beta\} = \{e, \alpha^n, \beta^3, \beta\},$$

**Proposition 4.3**

(1) $\langle \alpha^i, \gamma \rangle$ is a cyclic subgroups isomorphic to $C_{\frac{np}{i}}$ for all $i | 2n$

(2) $\langle \alpha^i, \alpha^j\beta, \gamma \rangle$ is a cyclic subgroups isomorphic to $C_{2p}$. where $i = n$ and $1 \leq j \leq n$.

*Proof.* A cyclic subgroup is clearly the type of subgroup that number (1) belongs to. It is still necessary to demonstrate that type (2)'s subgroup is cyclic.

$$(\alpha^j\beta\gamma)^{2p} = \underbrace{\alpha^j\beta\gamma . \alpha^j\beta\gamma \cdots \alpha^j\beta\gamma}_{2p \ times} = ((\alpha^j\beta)^2)^p((\gamma)^p)^2 = e$$

**Definition 4.5**

Let $G$ be a finite group. The following are hold:

1. $\Omega_G(\alpha) = \{\beta \in G | \langle \alpha, \beta \rangle \ is \ cyclic\}$, The cyclicizer of an element $x$ of $G$.
2. $\Omega_G(G) = \bigcap_{\alpha \in G} \Omega_G(\alpha) = \{\beta \in G | \langle \alpha, \beta \rangle \ is \ cyclic \ \forall \ \alpha \in G\}$

**Proposition 4.6**

Let $G$ be a finite group. The cyclicizer of an element $\alpha$ of $G$ is a subgroup of $G$ if and only if $G$ is a cyclic group.

**Proof:**

($\Rightarrow$) Assume the cyclicizer of an element $\alpha$ of $G$ is a subgroup of $G$ for each $\alpha$ in $G$

$\Omega_G(\alpha) = \{\beta \in G | \langle \alpha, \beta \rangle \ is \ cyclic\}$, thus mean $G$ is a cyclic group.

($\Longleftarrow$) Let $G$ be a cyclic group, then by definition the cyclicizer of any elements is define by $\Omega_G(\alpha) = \{\beta \in G | \langle \alpha, \beta \rangle \text{ is cyclic}\}$ is a subgroup of $G$.

**Theorem 4.7**

Let $n = 2^r \prod_{i \in \mathbb{N}} p_i^{\zeta_i}$ be integer number. The number of cyclic subgroups of group $\Phi$ is given by:

$$\Sigma CyS(\Phi) = \tau(p)(\tau(2n) + n)$$

**Proof:**

we first show that the set of all cyclic subgroups for each type of subgroups is presented in **Proposition 2.2**,

its clear that for each subgroup of $\langle \alpha^i \rangle$ and $i \mid 2n$ are equal to $\tau(2n)$, so for the subgroup of $\langle \alpha^i, c \rangle$ is same.

Now, the subgroup $\langle \alpha^i, \alpha^j \beta \rangle$ is a cyclic subgroup and isomorphic to $C_4$ if $i = n$, so so for the subgroup of $\langle \alpha^i, \alpha^j \beta, c \rangle$ is same.

**Example 4.8**

In the following example, we will present an application for computing the number of subgroups order table, using the Theorem [3.3]

Take $n = 2 3^2 5^2$ and $p = 7$, thus, the direct product group $\Phi = T_{1800} \times C_7$ thave order 12600, the number of all subgroups are given by the following table:

1) $k = 1,2,3,4$
2) The set of divisors of $n$ are $\{1,2,3,5,6,9,10,15,18,25,30,45,50,75,90,150,225,450\}$
3) Order of subgroups table:
4) If $\lambda|225$, then :

|  | k | 1 | 2 | 3 | 4 |
|---|---|---|---|---|---|
|  | * | 1 | 2 | $2^2$ | $2^3$ |
| $\lambda\|225$ |  | 1 | 1 | 1 | 451 | — |

|  | 3 | 1 | 1 | 1 | – |
|  | 5 | 1 | 1 | 1 | – |
|  | 9 | 1 | 1 | 1 | – |
|  | 15 | 1 | 1 | 1 | – |
|  | 25 | 1 | 1 | 1 | – |
|  | 45 | 1 | 1 | 1 | – |
|  | 75 | 1 | 1 | 1 | – |
|  | 225 | 1 | 1 | 1 | – |
| Sum First part | | 9 | 9 | 450 + 9 | |

5) If 7( $\lambda$|225), then :

| | k | 1 | 2 | 3 | 4 |
|---|---|---|---|---|---|
| | * | 1 | 2 | $2^2$ | $2^3$ |
| 7($\lambda$|225) | 7 | 1 | 1 | 451 | – |
| | 21 | 1 | 1 | 1 | – |

|   |   |   |   |   |
|---|---|---|---|---|
| 35 | 1 | 1 | 1 | — |
| 63 | 1 | 1 | 1 | — |
| 105 | 1 | 1 | 1 | — |
| 175 | 1 | 1 | 1 | — |
| 315 | 1 | 1 | 1 | — |
| 525 | 1 | 1 | 1 | — |
| 1575 | 1 | 1 | 1 | — |
| Sum Second Part | 9 | 9 | 450 + 9 | |

To clarify how these subgroups are calculated. We will calculate all subgroups based on the following Theorem[3.1]

From this type subgroup $\mathcal{H}_1 = \langle \alpha^i \rangle$ $\forall\, i|900$ we get the following subgroups:

| Subgroups | Order | Subgroups | Order | Subgroups | Order |
|---|---|---|---|---|---|
| $\langle \alpha \rangle$ | 900 | $\langle \alpha^{15} \rangle$ | 60 | $\langle \alpha^{75} \rangle$ | 12 |
| $\langle \alpha^2 \rangle$ | 450 | $\langle \alpha^{18} \rangle$ | 50 | $\langle \alpha^{90} \rangle$ | 10 |
| $\langle \alpha^3 \rangle$ | 300 | $\langle \alpha^{20} \rangle$ | 45 | $\langle \alpha^{100} \rangle$ | 9 |
| $\langle \alpha^4 \rangle$ | 225 | $\langle \alpha^{25} \rangle$ | 36 | $\langle \alpha^{150} \rangle$ | 6 |
| $\langle \alpha^5 \rangle$ | 180 | $\langle \alpha^{30} \rangle$ | 30 | $\langle \alpha^{180} \rangle$ | 5 |
| $\langle \alpha^6 \rangle$ | 150 | $\langle \alpha^{36} \rangle$ | 25 | $\langle \alpha^{225} \rangle$ | 4 |
| $\langle \alpha^9 \rangle$ | 100 | $\langle \alpha^{45} \rangle$ | 20 | $\langle \alpha^{300} \rangle$ | 3 |
| $\langle \alpha^{10} \rangle$ | 90 | $\langle \alpha^{50} \rangle$ | 18 | $\langle \alpha^{450} \rangle$ | 2 |
| $\langle \alpha^{12} \rangle$ | 75 | $\langle \alpha^{60} \rangle$ | 15 | $\langle e \rangle$ | 1 |

From this type subgroup $\langle \alpha^i, \alpha^j \beta \rangle$, $\forall\, i|450$ and $1 \leq j \leq i$;

we get the following subgroups:

$\langle \alpha^{450}, \alpha^j\beta \rangle_{1 \leq j \leq 450}$ are cyclic.

From this type subgroup $\mathcal{H}_3 = \langle \alpha^i, \gamma \rangle$ $\forall\, i|900$ we get the following subgroups:

| Subgroups | Order | Subgroups | Order | Subgroups | Order |
|---|---|---|---|---|---|
| $\langle \alpha, \gamma \rangle$ | 6300 | $\langle \alpha^{15}, \gamma \rangle$ | 420 | $\langle \alpha^{75}, \gamma \rangle$ | 84 |
| $\langle \alpha^2, \gamma \rangle$ | 3150 | $\langle \alpha^{18}, \gamma \rangle$ | 350 | $\langle \alpha^{90}, \gamma \rangle$ | 70 |
| $\langle \alpha^3, \gamma \rangle$ | 2100 | $\langle \alpha^{20}, \gamma \rangle$ | 315 | $\langle \alpha^{100}, \gamma \rangle$ | 63 |
| $\langle \alpha^4, \gamma \rangle$ | 1575 | $\langle \alpha^{25}, \gamma \rangle$ | 252 | $\langle \alpha^{150}, \gamma \rangle$ | 42 |
| $\langle \alpha^5, \gamma \rangle$ | 1260 | $\langle \alpha^{30}, \gamma \rangle$ | 210 | $\langle \alpha^{180}, \gamma \rangle$ | 35 |
| $\langle \alpha^6, \gamma \rangle$ | 1050 | $\langle \alpha^{36}, \gamma \rangle$ | 175 | $\langle \alpha^{225}, \gamma \rangle$ | 28 |
| $\langle \alpha^9, \gamma \rangle$ | 700 | $\langle \alpha^{45}, \gamma \rangle$ | 140 | $\langle \alpha^{300}, \gamma \rangle$ | 21 |
| $\langle \alpha^{10}, \gamma \rangle$ | 630 | $\langle \alpha^{50}, \gamma \rangle$ | 126 | $\langle \alpha^{450}, \gamma \rangle$ | 14 |
| $\langle \alpha^{12}, \gamma \rangle$ | 525 | $\langle \alpha^{60}, \gamma \rangle$ | 105 | $\langle \gamma \rangle$ | 7 |

From this type subgroup $\langle \alpha^i, \alpha^j\beta, \gamma \rangle$, $\forall\, i|450$ and $1 \leq j \leq i$;

we get the following subgroups:

$$\langle \alpha^{450}, \alpha^j\beta, \gamma \rangle_{1 \leq j \leq 450}$$

Sum first part = 1236 subgroups

Sum second part = 1236 subgroups

So, we found out that there are a total of 1236+1236=2472 subgroups in $\Phi = T_{1800} \times C_7$.

By using the final formula proposed in the Theorem 3.2 to verify the final result, we notice that the number be equal to the number of calculation in the table.

$$\Sigma NS(T_{1800} \times C_7) = \tau(7)(\tau(900) + \sigma(450))$$
$$= (2)(27 + 1209)$$
$$= 2472$$

By GAP program.

```
gap> n:=450;

gap>F := FreeGroup( "a", "b" );

gap>T := F / [ a^(2*n), b^2/a^n, b^(-1)*a*b*a ];

gap> i:=IsomorphismPermGroup(T);

gap>c:=CyclicGroup(7);

gap> Tc:=DirectProduct(T,c);

gap> s:=AllSubgroups(Tc);

gap>x:=[];

gap> for t in s do

gap> dd:=IsCyclic(t);

gap> if dd=true then; Add(x,t), fi; od; x; Size(x);
```

Answer: 2472

## 5- Aplication.

In our work, suppose that the our group is define by:

$T_{4n} \times C_p = = \langle \alpha, \beta, \gamma | \alpha^{2n} = e, \alpha^n = \beta^2 = \gamma^p = e, \beta\alpha\beta^{-1} = \alpha^{-1}, [\alpha, \gamma] = [\beta, \gamma] = e \rangle$

be a finite group has order $4np$ and we denoted by $\Sigma CyS(G)$ is a set of all cyclic subgroups of $T_{4n} \times C_p$. Here, we work on an important field of group theory. Tărnăuceanu [7,18,19] proved that a finite group $G$ has $|G| - 1$ cyclic subgroups if and only if $G$ is isomorphic to $C_k, k = 3,4$ or it isomorphic to $D_{2k}, k = 3,4$. In this paper, Describe the finite groups $T_{4n} \times C_p$ satisfying $\Sigma CyS(G) = |T_{4n} \times C_p| - t$.

**Proposition 5.1**

Suppose that where $p$ and $q$ are distinct primes and $p < q$ and $T_{4n} \times C_p$ be a finite non-abelian group of order $4qp$. If $n = q$, then $t = 2q(2p - 1) - 8$

$$\Sigma CuS(T_{4n} \times C_p) = 4pq - 2q(2p-1) - 8$$

**Proof:**

Let $n = q$ be a prime number, and p is prime number and $p < q$,

By applying Theorem 4.8, we obtain on the following.

$$|T_{4n} \times C_p| = 4qp$$

$$\Sigma CyS(T_{4n} \times C_p) = \tau(p)(\tau(2q) + q)$$

$$\Sigma CyS(T_{4n} \times C_p) = 2(\tau(2)\tau(q) + q)$$

$$\Sigma CyS(T_{4n} \times C_p) = 2(4 + q)$$

$$t = 4qp - 2(4 + q)$$

$$t = 2q(2p - 1) - 8$$

The proof is done.

**Proposition 5.2**

Suppose that where $p$ and $q^2$ are distinct primes and $p < q$ and $T_{4n} \times C_p$ be a finite non-abelian group of order $4q^2p$. If $n = q^2$, then $t = 2q^2(2p-1) - 12$

$$\Sigma CyS(T_{4n} \times C_p) = 4q^2p - 2q^2(2p-1) - 12$$

**Proof:**

Let $n = q^2$ be a prime number, and p is prime number and $p < q$,

By applying Theorem 4.8, we obtain on the following.

$$|T_{4n} \times C_p| = 4q^2p$$

$$\Sigma CyS(T_{4n} \times C_p) = \tau(p)(\tau(2q^2) + q^2)$$

$$\Sigma CyS(T_{4n} \times C_p) = 2(\tau(2)\tau(q^2) + q^2)$$

$$\Sigma CyS(T_{4n} \times C_p) = 2(6 + q^2)$$

$$t = 4q^2p - 2(6 + q^2)$$

$$t = 2q^2(2p - 1) - 12$$

The proof is done.

**Theorem 5.3**

Let $n = q^r$ be a prime number, $r \geq 1$ and $p < q$, then

$$t = 2q^r(2p - 1) - 4\tau(q^r)$$

**Proof.**

Let $n = q^r$ be a prime number, and p is prime number and $p < q$,

By applying Theorem 4.8, we obtain on the following.

$$|T_{4n} \times C_p| = 4q^r p$$

$$\Sigma CyS(T_{4n} \times C_p) = \tau(p)(\tau(2q^r) + q^r)$$

$$\Sigma CyS(T_{4n} \times C_p) = 2(\tau(2)\tau(q^r) + q^r)$$

$$\Sigma CyS(T_{4n} \times C_p) = 2(2(r + 1) + q^r)$$

$$t = 4q^r p - 2(2(r + 1) + q^2)$$

$$t = 2q^r(2p - 1) - 2(2(r + 1))$$

$$t = 2q^r(2p - 1) - (4(r + 1))$$

$$t = 2q^r(2p - 1) - (4\tau(q^r))$$

**Theorem 5.4**

Let $n = 2^r$ be an even prime number, $r \geq 1$ and $p < q$, then

$$t = 2^r(2p - 1) - 2\tau(2^{r+1})$$

**Proof.**

Let $n = 2^r$ be an even prime number, and p is prime number and $p < q$,

By applying Theorem 4.8, we obtain on the following.

$$|T_{4n} \times C_p| = 42^r p = 2^{r+2} p$$

$$\Sigma CyS(T_{4n} \times C_p) = \tau(p)(\tau(22^r) + 2^r)$$

$$\Sigma CyS(T_{4n} \times C_p) = 2(\tau(2^{r+1}) + 2^r)$$

$$\Sigma CyS(T_{4n} \times C_p) = 2((r + 2) + 2^r)$$

$$t = 2^{r+2} p - 2((r + 2) + 2^r)$$

$$t = 2^r(2^2 p - 1) - 2((r + 2))$$

$$t = 2^r(2^2 p - 1) - (2\tau(2^{r+1}))$$

**Theorem 5.5**

If $n$ be a positive integer number, then

$$\Sigma NS(T_{4n} \times C_p) \leq \Sigma CyS(T_{4n} \times C_p)$$

**Proof.**

In this proof, we must spilt to two cases:

Case 1, Suppose that $n = 1$, then it is true.

$$\Sigma NS(T_{4n} \times C_p) = \Sigma CyS(T_{4n} \times C_p)$$

Case 2,

Suppose that $n > 1$, thus,

$$\Sigma NS(T_{4n} \times C_p) = \tau(p)(\tau(2n) + \sigma(n))$$

and

$$\Sigma CyS(T_{4n} \times C_p) = \tau(p)(\tau(2n) + n$$

Suppose that are equal, we obtain on the

$$\tau(p)(\tau(2n) + \sigma(n)) = \tau(p)(\tau(2n) + n)$$

By using cancel law obtain on the $\sigma(n) = n$ it is not true for $n > 1$,

Thus $\sigma(n) < n$.

## Corollary 5.6

If $n$ be a positive integer number, then

$$\Sigma NS(T_{4n} \times C_p) \leq |T_{4n} \times C_p| - t$$

For all $t \in \mathbb{N}$.

**Proof:**

Directly from Theorem 5.5

$$\Sigma NS(T_{4n} \times C_p) \leq \Sigma CyS(T_{4n} \times C_p)$$

and by theorem 2.1 in [16]

$$\Sigma CyS(G) = |T_{4n} \times C_p| - t$$

We obtain on the

$$\Sigma NS(T_{4n} \times C_p) \leq \Sigma CyS(T_{4n} \times C_p) = |T_{4n} \times C_p| - t$$

$$\Sigma NS(T_{4n} \times C_p) \leq |T_{4n} \times C_p| - t$$